\documentstyle[12pt]{amsart}
\begin{document}
\title{The biharmonicity of sections of the tangent bundle}
\author{M.~Markellos and H.~Urakawa}
\title[The biharmonicity]
{The biharmonicity of sections of the tangent bundle}


\dedicatory{}
\maketitle
\begin{abstract}
The bienergy of a vector field on a Riemannian manifold $(M, g)$ is defined to be the bienergy of the corresponding map $(M, g)\mapsto (TM, g_{S})$, where the tangent bundle $TM$ is equipped with the Sasaki metric $g_{S}$. The constrained variational problem is studied, where variations are confined to vector fields, and the corresponding critical point condition characterizes biharmonic vector fields. Furthermore, we prove that if $(M, g)$ is a compact oriented $m$-dimensional Riemannian manifold and $X$ a tangent vector of $M$, then $X$ is a biharmonic vector field of $(M, g)$ if and only if $X$ is parallel. Finally, we give examples of non-parallel biharmonic vector fields in the case which the base manifold $(M, g)$ is non-compact.
\end{abstract}

\numberwithin{equation}{section}
\theoremstyle{plain}
\newtheorem{df}{Definition}[section]
\newtheorem{th}[df]{Theorem}
\newtheorem{prop}[df]{Proposition}
\newtheorem{lem}[df]{Lemma}
\newtheorem{cor}[df]{Corollary}
\newtheorem{rem}[df]{Remark}
\newtheorem{exam}[df]{Example}

\section{Introduction}
The energy functional of a map $\varphi:(M, g)\mapsto (N, h)$ between
Riemannian manifolds has been widely investigated by several
researchers (\cite{Baird03}, \cite{EellsSampson64}, \cite{Urakawa})
and is given by
\begin{equation*}
E_{1}(\varphi)=\frac{1}{2}\int_{M}{\|d\varphi\|^{2}v_{g}}
\end{equation*}
where $d\varphi$ denotes the differential of the map $\varphi$. Critical
points for the energy functional are called \textit{harmonic maps}
and have been characterized by the vanishing of the \textit{tension
field} $\tau_{1}(\varphi)={\rm Tr} \nabla d\varphi$.

Let $(M, g)$ be a Riemannian manifold and denote by $(TM,
g_{S})$ its tangent bundle equipped with the Sasaki
metric $g_{S}$. A vector field $X$ on $M$ determines a mapping
from $(M, g)$ into $(TM, g_{S})$, embedding $M$ into its
tangent bundle. If $M$ is compact and orientable, the energy of $X$ is the energy of the corresponding map (\cite[pp. 41-42]{DragomirPerro11}). Nouhaud (\cite{Nouhaud77}) proved that if $M$ is compact, then $X: (M, g)\mapsto (TM, g_{S})$ is an harmonic map if and only if $X$ is parallel. The same conclusion is obtained in \cite{Ishiha79} but, in this case, the explicit expression of the tension field of the vector field is provided (see subsection 2.2). Gil - Medrano (\cite{GilMedr01}) proved that if $M$ is compact, then a vector field $X$ is a critical point of the energy functional $E_{1}$ restricted to the set ${\mathfrak{X}}(M)$ of all vector fields of $(M, g)$ (equivalently, $X$ is a \emph{harmonic vector field}) if and only if $X$ is parallel. As a consequence, the search for critical points for the energy functional $E_{1}:C^{\infty}(M, TM)\mapsto [0, +\infty)$ or $E_{1}:{\mathfrak{X}}(M)\mapsto [0, +\infty)$ shows that both domain $C^{\infty}(M, TM)$ and ${\mathfrak{X}}(M)$ are inappropriate.

In \cite{EellsSampson64}, J.Eells and J.H.Sampson suggested the idea
of studying $k$-harmonic maps. A particular interest has the case
$k=2$. They define the \textit{bienergy} of $\varphi$ as the functional
\begin{equation*}
E_{2}(\varphi)=\frac{1}{2}\int_{M}{\|\tau_{1}(\varphi)\|^{2}v_{g}},
\end{equation*}
and a map is \textit{biharmonic} if and only if it is a critical
point of $E_{2}$. Jiang (\cite{Jiang08}) derived the associated
Euler-Lagrange equation for $E_{2}$. A harmonic map is automatically a biharmonic map. Non-harmonic biharmonic maps are said to be
\textit{proper biharmonic} maps. In the last decade there has been a growing interest in the theory of biharmonic maps which can be divided in two main research directions. On the one side, the differential geometric aspect has driven attention to the construction of examples and classification results. More precisely,  in \cite{CaMoOn01}, Caddeo et al. classified biharmonic curves and surfaces of the unit 3-sphere ${\mathbb{S}}^{3}$. In fact, they found that they
are circles, helices which are geodesics in the Clifford minimal torus and small hyperspheres. The same authors in \cite{CaMoOn02}
constructed examples of proper biharmonic submanifolds of ${\mathbb{S}}^{n}$, $n>3$. The other side is the analytic aspect from the point view of PDE: biharmonic maps are solutions of a fourth order strongly elliptic semilinear PDE.

We concentrate on the mapping $X:(M, g)\mapsto (TM, g_{S})$. It is natural to consider the problem of characterizing those vector fields for which the corresponding map is a biharmonic map. Furthermore, we can also look for vector fields $X$ that are critical points of the bienergy functional restricted to variations among vector fields. In the latter case, such a vector field is called \emph{biharmonic vector field}. The goal of this paper is to answer to these problems. The paper is organized in the following way. Section 2 contains the presentation of some basic notions about the geometry of the tangent bundle and biharmonic maps. In Section 3, we derive the first variational formula associated to the bienergy functional restricted to the space of all vector fields of $(M, g)$ (see Theorem 3.2). By using this formula, we prove the following Theorem:
\begin{th}
Let $(M, g)$ be a compact $m$-dimensional Riemannian manifold and $X\in{\mathfrak{X}}(M)$ a tangent vector field. Then, $X$ is a biharmonic vector field if and only if $X$ is parallel.
\end{th}
As a consequence, we deduce that
\begin{cor}
Let $(M, g)$ be a compact oriented $m$-dimensional Riemannian manifold and $X\in{\mathfrak{X}}(M)$ a tangent vector field. Then, $X:(M, g)\mapsto (TM, g_{S})$ is a biharmonic map if and only if $X$ is parallel.
\end{cor}
Finally, we give examples of non-parallel biharmonic vector fields of $({\mathbb{R}}^{2}, g)$, equipped with the standard Euclidean metric $g$, which are non harmonic (see Example 3.9 and Remarks 3.10 and 3.11). More precisely, in Remark 3.11, we give the general representation formula of a biharmonic vector field of $({\mathbb{R}}^{2}, g)$. Furthermore, we give examples of biharmonic vector fields which are non harmonic vector fields in the case which the base manifold $(M, g)$ is non-flat (see Examples 3.12 and 3.13). As a consequence, we point out that the notions ''biharmonic vector field'' and ''harmonic vector field'' are independent in the sense that it does not exist an immediate relation between them (\footnote{In \cite{DjaaElheOuak12}, the authors proved Corollary 1.2 following a quite different approach from that of the paper. We should point out that their paper includes some serious mistakes. More precisely, on page 470 and line 14, the term $\nabla^{X}_{e_{i}}(\tau^{v}(X))^{V}$should be calculated at the point $(x, X_{x})\in TM$ in order to be a section of the pull - back bundle $X^{-1}(TTM)$. As a consequence, this mistake is transferred to the calculation of the term ${\rm Tr}_{g}\nabla^{2}(\tau^{v}(X))^{V}_{(x, u)}$ (page 470 and line 18)}).


\section{Preliminaries.}
\subsection{The tangent bundle} We recall some basic facts about the geometry of the tangent bundle.
For a more elaborate exposition, we refer to the survey \cite{GudmuKappos02}. In the present paper, we denote by
${\mathfrak{X}}(M)$ the space of all vector fields of a Riemannian manifold $(M, g)$.
\par
Let $(M, g)$ be an $m$-dimensional Riemannian manifold and $\nabla$
the associated Levi-Civita connection. Its Riemann curvature tensor
$R$ is defined by
\begin{equation*}
R(X, Y)Z =\nabla_{X}\nabla_{Y}Z - \nabla_{Y}\nabla_{X}Z -
\nabla_{[X, Y]}Z
\end{equation*}
for all vector fields $X, Y$ and $Z$ on $M$. The tangent bundle of a
Riemannian manifold $(M, g)$ is denoted by $TM$ and consists of
pairs $(x, u)$ where $x$ is a point in $M$ and $u$ a tangent vector
to $M$ at $x$. The mapping $\pi: TM\mapsto M: (x, u)\mapsto x$ is
the natural projection from $TM$ onto $M$. The tangent space $T_{(x, u)}TM$ at a point $(x,
u)$ in $TM$ is a direct sum of the vertical subspace ${\mathcal{V}}_{(x,
u)}={\rm Ker}(d\pi|_{(x, u)})$ and the horizontal subspace ${\mathcal{H}}_{(x,
u)}$, with respect to the Levi-Civita connection $\nabla$ of $M$:
\begin{equation*}
T_{(x, u)}TM={\mathcal{H}}_{(x, u)}\oplus {\mathcal{V}}_{(x, u)}.
\end{equation*}
For any vector $w\in T_{x}M$, there exists a unique vector $w^{h}\in
{\mathcal{H}}_{(x, u)}$ at the point $(x, u)\in TM$, which is called the
\emph{horizontal lift} of $w$ to $(x, u)$, such that $d\pi(w^{h})=w$
and a unique vector $w^{v}\in {\mathcal{V}}_{(x, u)}$, which is called the
\emph{vertical lift} of $w$ to $(x, u)$, such that $w^{v}(df)=w(f)$
for all functions $f$ on $M$. Hence, every tangent vector
$\bar{w}\in T_{(x, u)}TM$ can be decomposed as
$\bar{w}=w_{1}^{h}+w_{2}^{v}$ for uniquely determined vectors
$w_{1}, w_{2} \in T_{x}M$. The \emph{horizontal} (respectively,
\emph{vertical}) \emph{lift} of a vector field $X$ on $M$ to $TM$ is
the vector field $X^{h}$ (respectively, $X^{v}$) on $TM$ whose value
at the point $(x, u)$ is the horizontal (respectively, vertical)
lift of $X_{x}$ to $(x, u)$.

The tangent bundle $TM$ of a Riemannian manifold $(M, g)$ can be
endowed in a natural way with a Riemannian metric $g_{S}$, the
\emph{Sasaki metric}, depending only on the Riemannian structure $g$
of the base manifold $M$. It is uniquely determined by
\begin{equation}\label{Eq: Sasaki metric}
\begin{array}{ll}
g_{S}(X^{h}, Y^{h})=g_{S}(X^{v}, Y^{v})=g(X, Y)\circ \pi, &g_{S}(X^{h}, Y^{v})=0\\
\end{array}
\end{equation}
for all vector fields $X$ and $Y$ on $M$. More intuitively, the
metric $g_{S}$ is constructed in such a way that the vertical and
horizontal subbundles are orthogonal and the bundle map $\pi: (TM,
g_{S})\mapsto (M, g)$ is a Riemannian submersion.
\par

\subsection{Biharmonic maps}
Let $(M, g), (N, h)$ be Riemannian manifolds of dimensions $m$ and $n$, respectively, and let $\varphi:
(M, g)\mapsto (N, h)$ be a smooth map between them. We denote by $\nabla^{\varphi}$ the connection
of the vector bundle $\varphi^{-1}TN$ induced from the Levi-Civita
connection $\bar{\nabla}$ of $(N, h)$ and $\nabla$ the Levi-Civita
connection of $(M, g)$. Let $D$ be a compact domain of $M$. The \textit{energy (integral) of
$\varphi$} over $D$ is defined by:
\begin{equation*}
E_{1}(\varphi)=\frac{1}{2}\int_{D}\|d\varphi\|^{2}v_{g},
\end{equation*}
where $v_{g}$ is the volume element of $(M, g)$.

A smooth map $\varphi: (M, g)\mapsto (N, h)$ is said to be \textit{harmonic} if it is a critical point of the energy functional
for any compact domain $D$. It is well known (\cite{Baird03}) that the smooth map $\varphi:
(M, g)\mapsto (N, h)$ is harmonic if and only if
\begin{equation*}
\tau_{1}(\varphi)={\rm Tr}(\nabla
d\varphi)=\sum_{i=1}^{m}\{\nabla_{e_{i}}^{\varphi}d\varphi(e_{i})-d\varphi(\nabla_{e_{i}}{e_{i}})\}=0,
\end{equation*}
where $\{e_{i}\}_{i=1}^{m}$ is a local orthonormal frame field of $(M, g)$. The
equation $\tau_{1}(\varphi)=0$ is called the \emph{harmonic equation}.

J.~Eells and J.~H.~Sampson \cite{EellsSampson64} introduced the
notion of polyharmonic maps. In this paper, we only consider
polyharmonic maps of order two. Such maps are frequently called
\emph{biharmonic maps}.

More precisely, a smooth map $\varphi: (M, g)\mapsto
(N, h)$ is said to be \emph{biharmonic} if it is a critical
point of the bienergy functional:
\begin{equation*}
E_{2}(\varphi)=\frac{1}{2}\int_{D}\|\tau_{1}(\varphi)\|^{2}v_{g},
\end{equation*}
over every compact region $D$ of $M$. The corresponding Euler-Lagrange equation associated to the bienergy functional
becomes more complicated and it involves the curvature of the
target manifold $(N, h)$ (\cite{Jiang08}). We should mention that an harmonic map is automatically a
biharmonic map, in fact a minimum of the bienergy functional. We refer to the survey \cite{MontaOni06} for more information on results related to the theory of biharmonic maps.

A vector field $X$ on $(M, g)$ can be regarded as the immersion
$X: (M, g)\mapsto (TM, g_{S}): x\mapsto (x, X_{x})\in TM$
into its tangent bundle $TM$ equipped with the
Sasaki metric $g_{S}$. The tension field $\tau_{1}(X)$ is given by (\cite{GilMedr01})
\begin{equation}\label{Eq: tension field of X}
\tau_{1}(X)=\Big(-\sum_{i=1}^{m}{R(\nabla_{e_{i}}{X}, X)e_{i}}\Big)^{h}
+({\rm Tr} \nabla^{2}X)^{v},
\end{equation}
where $\{e_{i}\}_{i=1}^{m}$ is a local orthonormal frame field of $(M, g)$. We
mention that the horizontal part of $\tau_{1}(X)$ is appeared with
different sign from the same expression in (\cite{GilMedr01}) due to
the different sign of the Riemann curvature tensor. The term
$-{\rm Tr}\nabla^{2}X$ is equal to the \emph{rough Laplacian}
$\bar{\Delta}X$, that is
\begin{equation*}
\bar{\Delta}X=-{\rm Tr}\nabla^{2}X=\sum_{i=1}^{m}\{\nabla_{\nabla_{e_{i}}{e_{i}}}X-\nabla_{e_{i}}\nabla_{e_{i}}X\}.
\end{equation*}
For the sake of convenience, we set
$S(X)=\sum_{i=1}^{m}{R(\nabla_{e_{i}}{X}, X)e_{i}}$. It is worth noting that the term $S(X)$ is
well-defined independently of the orthonormal frame used. Furthermore, it holds
\begin{equation}\label{Eq: tensorial character S(X)}
S(fX)=f^{2}S(X),
\end{equation}
for any smooth function $f$ and vector field $X$ of $M$. As a consequence, relation (\ref{Eq: tension field of X}) is transformed to
\begin{equation}\label{Eq: final tension of X}
\tau_{1}(X)=(-S(X))^{h}+(-\bar{\Delta}X)^{v}.
\end{equation}
In the following, we will use the following Lemmas:
\begin{lem}
Let $(M, g)$ be an $m$-dimensional Riemannian manifold and $X, Y$ vector fields on $M$. Then the following equation is satisfied:
\begin{equation}\label{Eq: helpful lemma}
\Delta(g(X, Y))=g(\bar{\Delta}X, Y)-g(X, \bar{\Delta}Y)-2div\theta_{XY},
\end{equation}
where $\theta_{XY}$ is the 1-form defined by $\theta_{XY}(Z)=g(X, \nabla_{Z}Y)$, ($Z\in {\mathfrak{X}}(M)$) and $\Delta$ is the ordinary Laplace-Beltrami operator acting on functions.
\end{lem}
\begin{pf}
Let $\{e_{i}:1\leq i \leq m\}$ be a local orthonormal frame field of $(M, g)$. Then (due to the definition of the rough Laplacian and $\nabla g=0$)
\begin{eqnarray}
\Delta(g(X, Y))&=&\sum_{i=1}^{m}\Big\{\nabla_{e_{i}}e_{i}(g(X, Y))-e_{i}e_{i}(g(X, Y))\Big\}\nonumber\\
&=&\sum_{i=1}^{m}\Big\{g(\nabla_{\nabla_{e_{i}}e_{i}}X, Y)-g(\nabla_{e_{i}}\nabla_{e_{i}}X, Y)+g(X, \nabla_{\nabla_{e_{i}}e_{i}}Y)-\nonumber\\
&&-g(X, \nabla_{e_{i}}\nabla_{e_{i}}Y)-2g(\nabla_{e_{i}}X, \nabla_{e_{i}}Y)\Big\}=g(\bar{\Delta}X, Y)+\nonumber\\
&&+g(X, \bar{\Delta}Y)-2\sum_{i=1}^{m}g(\nabla_{e_{i}}X, \nabla_{e_{i}}Y)=g(\bar{\Delta}X, Y)+\nonumber\\
&&+g(X, \bar{\Delta}Y)-2\sum_{i=1}^{m}\Big\{e_{i}(g(X, \nabla_{e_{i}}Y))-g(X, \nabla_{\nabla_{e_{i}}e_{i}}Y)+\nonumber\\
&&+g(X, \nabla_{\nabla_{e_{i}}e_{i}}Y)-g(X, \nabla_{e_{i}}\nabla_{e_{i}}Y)\Big\}=g(\bar{\Delta}X, Y)+\nonumber\\
&&+g(X, \bar{\Delta}Y)-2div\theta_{XY}-2g(X, \bar{\Delta}Y)=g(\bar{\Delta}X, Y)-\nonumber\\
&&-g(X, \bar{\Delta}Y)-2div\theta_{XY}.\nonumber
\end{eqnarray}
Here we made use of the local expression (with respect to the local orthonormal frame field $\{e_{i}:1\leq i \leq m\}$ of the divergence of differential 1-forms i.e.,
\begin{equation*}
\begin{array}{lr}
div\omega=\sum_{i=1}^{m}(\nabla_{e_{i}}\omega)e_{i}, & \omega \in \Omega^{1}(M).
\end{array}
\end{equation*}
\end{pf}
\begin{rem}
Lemma 2.1 could be considered as a generalization of the Lemma 2.15 of \cite[p. 53]{DragomirPerro11}.
\end{rem}
\begin{lem}\cite{HiKaWo01}
Let $(M, g)$ be a Riemannian manifold and $X$ a vector field of M.
Then the following equation is satisfied:
\begin{equation}\label{Eq: rough Laplacian fX}
\bar{\Delta}(fX)=(\Delta f) X +f \bar{\Delta}X-2\nabla_{grad f}X,
\end{equation}
where $f$ being a smooth function of $M$ and $grad f$ the gradient of
$f$.
\end{lem}

\section{The biharmonicity of vector fields}

In this section, we define the bienergy $E_{2}(X)$ of a vector field $X$ of a Riemannian manifold $(M, g)$.
Furthermore, we determine the first variational formula of the bienergy functional $E_{2}$ restricted to
the space of all vector fields.

Let $(M, g)$ be a compact $m$-dimensional Riemannian manifold, $(TM, g_{S})$ its tangent bundle
equipped with the Sasaki metric $g_{S}$ and $X$ a tangent vector field on $M$. The bienergy $E_{2}(X)$ of
$X$ is defined to be the bienergy of the corresponding map $X:(M, g)\mapsto (TM, g_{S})$. More precisely,
combining relations (\ref{Eq: Sasaki metric}) and (\ref{Eq: final tension of X}), we get
\begin{equation}\label{Eq: bienergy functional}
E_{2}(X)=\frac{1}{2}\int_{M}\|\tau_{1}(X)\|^{2}v_{g}=\frac{1}{2}\int_{M}\Big[g(S(X), S(X))+g(\bar{\Delta}X, \bar{\Delta}X)\Big]v_{g}.
\end{equation}
Now, we give the following definition:
\begin{df}
Let $(M, g)$ be a Riemannian manifold. A vector field $X\in{\mathfrak{X}}(M)$ is called \emph{biharmonic} if the corresponding map $X:(M, g)\mapsto (TM, g_{S})$ is a critical point for the bienergy functional $E_{2}$, only considering variations among maps defined by vector fields.
\end{df}
In the sequel, we deduce the critical point condition for $E_{2}:{\mathfrak{X}}(M)\mapsto {\mathbb{R}}$. We prove the following Theorem:
\begin{th}
Let $(M, g)$ be a compact oriented $m$-dimensional Riemannian manifold, $\{e_{i}\}_{i=1}^{m}$ a local orthonormal frame field of $(M, g)$, $X$ a tangent vector field on $M$ and $E_{2}:{\mathfrak{X}}(M)\mapsto [0, +\infty)$ the bienergy functional restricted to the space of all vector fields. Then
\begin{eqnarray}\label{Eq: variational formula}
\frac{d}{dt}E_{2}(X_{t})_{t=0}&=&\int_{M}\Big\{g(\bar{\Delta}\bar{\Delta}X+\sum_{i=1}^{m}\Large[(\nabla_{e_{i}}R)(e_{i}, S(X))X+R(e_{i},
\nabla_{e_{i}}S(X))X+\nonumber\\
&&+2R(e_{i}, S(X))\nabla_{e_{i}}X\Large], V)\Big\}v_{g}
\end{eqnarray}
for any smooth 1-parameter variation $U:M\times (-\epsilon, \epsilon)\mapsto TM$ of $X$ through vector fields i.e., $X_{t}(z)=U(z, t)\in T_{z}M$ for any $|t|<\epsilon$ and $z\in M$, or equivalently $X_{t}\in{\mathfrak{X}}(M)$ for any $|t|<\epsilon$. Also, $V$ is the tangent vector field on $M$ given by
\begin{equation*}
\begin{array}{lr}
V(z)=\frac{d}{dt}X_{z}(0), & z\in M,
\end{array}
\end{equation*}
where $X_{z}(t)=U(z, t), (z, t)\in M\times (-\epsilon, \epsilon)$.
\end{th}
\begin{pf}
Let $X$ be a tangent vector on $M$ and $I=(-\epsilon, \epsilon), \epsilon>0$. For $t\in I$,
we denote by $i_{t}:M\mapsto M\times I, p\mapsto (p,t)$ the
canonical injection. We consider $C^{\infty}$-variations $U:
M\times I\mapsto TM$ of $X$ within ${\mathfrak{X}}(M)$, i.e. for all
$t\in I$ the mappings $X_{t}=U\circ i_{t}$ are in fact vector
fields and $X_{0}=X$. We choose $\{e_{i}\}_{i=1}^{m}$ a local orthonormal
frame field of $(M, g)$. We extend $e_{i}$ (resp. $\frac{d}{dx}\in {\mathfrak{X}}(I)$) to $M\times I$, denoted by $E_{i}$(resp. $\frac{d}{dt}$).
Moreover, we have $[E_{i}, \frac{d}{dt}]=0$. We denote by $D$ the
Levi-Civita connection of $M\times I$ and $R^{D}$ the Riemann
curvature tensor of $M\times I$. Since $M\times I$ is a Riemannian product, we have (using the second Bianchi identity for the last relation)
\begin{equation}\label{Eq: Riemannian product}
\begin{array}{llcr}
R^{D}(TM, TI)=0, &D_{\frac{d}{dt}}E_{i}=0,  & D_{E_{i}}\frac{d}{dt}=0, &(D_{\frac{d}{dt}}R^{D})(D_{E_{i}}U, U)E_{i}=0
\end{array}
\end{equation}
for all $1\leq i \leq m$. We set $Z=\sum_{i=1}^{m} R^{D}(D_{E_{i}}U, U)E_{i}$ and
$\Omega=\sum_{i=1}^{m}[D_{D_{E_{i}}E_{i}}U-D_{E_{i}}D_{E_{i}}U]$. We easily observe that $S(X_{t})=Z\circ i_{t}$ and $\bar{\Delta}X_{t}=\Omega\circ i_{t}$. In the sequel, we consider the function
\begin{eqnarray}
E_{2}(t)&=&E_{2}(X_{t})=\frac{1}{2}\int_{M}[g(S(X_{t}), S(X_{t}))+g(\bar{\Delta}X_{t}, \bar{\Delta}X_{t})]v_{g}\nonumber\\
&=&\frac{1}{2}\int_{M}[g(Z, Z)+g(\Omega, \Omega)]\circ i_{t}v_{g}.\nonumber
\end{eqnarray}
Differentiating the function $E_{2}(t)$ at each $t$, we obtain
\begin{equation}\label{Eq: derivative E2}
\frac{d}{dt}E_{2}(X_{t})=\int_{M}g(D_{\frac{d}{dt}}Z, Z)\circ
i_{t}v_{g}+\int_{M}g(D_{\frac{d}{dt}}\Omega, \Omega)\circ
i_{t}v_{g}.
\end{equation}
Taking into account the symmetries of the Riemann curvature tensor,
relations (\ref{Eq: Riemannian product}) and summing over all repeated indices, we have
\begin{eqnarray}\label{Eq: derivative of Z}
\int_{M}g(D_{\frac{d}{dt}}Z, Z)\circ
i_{t}v_{g}&=&\int_{M}g\Big ((D_{\frac{d}{dt}}R^{D})(D_{E_{i}}U,
U)E_{i}+R^{D}(D_{\frac{d}{dt}}D_{E_{i}}U,
U)E_{i}+\nonumber\\
&&+R^{D}(D_{E_{i}}U, D_{\frac{d}{dt}}U)E_{i}, Z\Big)\circ
i_{t}v_{g}\nonumber\\
&=&\int_{M}\Big[g\Big(R^{D}
(D_{E_{i}}D_{\frac{d}{dt}}U+R^{D}(\frac{d}{dt},
E_{i})U, U)E_{i}, Z\Big)+\nonumber\\
&&+g(R^{D}(E_{i}, Z)D_{E_{i}}U,
D_{\frac{d}{dt}}U)\Big]\circ i_{t}v_{g}\nonumber\\
&=&\int_{M}\Big[-g(R^{D}(E_{i}, Z)U,
D_{E_{i}}D_{\frac{d}{dt}}U)+\nonumber\\
&&+g(R^{D}(E_{i}, Z)D_{E_{i}}U,
D_{\frac{d}{dt}}U)\Big]\circ i_{t}v_{g}\nonumber\\
&=&\int_{M}\Big\{[-E_{i}(g(R^{D}(E_{i}, Z)U, D_{\frac{d}{dt}}U))+g(R^{D}(D_{E_{i}}E_{i}, Z)U, D_{\frac{d}{dt}}U)]+\nonumber\\
&&+g((D_{E_{i}}R^{D})(E_{i}, Z)U, D_{\frac{d}{dt}}U)+g(R^{D}(E_{i}, D_{E_{i}}Z)U, D_{\frac{d}{dt}}U)+\nonumber\\
&&+2g(R(E_{i}, Z)D_{E_{i}}U, D_{\frac{d}{dt}}U)\Big\}\circ i_{t}v_{g}.
\end{eqnarray}
Applying the divergence Theorem for the 1-form $\eta_{t}(W)=g(R(W, S(X_{t}))X_{t}, \nabla_{\frac{d}{dt}}X_{t}), t\in I,  W\in {\mathfrak{X}}(M)$, relation (\ref{Eq: derivative of Z}) gives
\begin{eqnarray}\label{Eq: final derivative of Z}
\int_{M}g(D_{\frac{d}{dt}}Z, Z)\circ i_{t}v_{g}&=&\int_{M}g\Big((\nabla_{e_{i}}R)(e_{i}, S(X_{t}))X_{t}+g(R(e_{i}, \nabla_{e_{i}}S(X_{t}))X_{t}+\nonumber\\
&&+2g(R(e_{i}, S(X_{t}))\nabla_{e_{i}}X_{t}, \nabla_{\frac{d}{dt}}X_{t}\Big)v_{g}.
\end{eqnarray}
Similarly, summing over all repeated indices, we deduce
\begin{eqnarray}
\int_{M}g(D_{\frac{d}{dt}}\Omega, \Omega)\circ
i_{t}v_{g}&=&\int_{M}g(D_{\frac{d}{dt}}D_{D_{E_{i}}E_{i}}U-D_{\frac{d}{dt}}D_{E_{i}}D_{E_{i}}U, \Omega)\circ i_{t}v_{g}\nonumber\\
&=&\int_{M}g(D_{D_{E_{i}}E_{i}}D_{\frac{d}{dt}}U-D_{E_{i}}D_{E_{i}}D_{\frac{d}{dt}}U, \Omega)\circ i_{t}v_{g}\nonumber\\
&=&\int_{M}\Big \{D_{E_{i}}E_{i}[g(D_{\frac{d}{dt}}U, \Omega)]-g(D_{\frac{d}{dt}}U, D_{D_{E_{i}}E_{i}}\Omega)-\nonumber\\
&&-E_{i}[g(D_{E_{i}}D_{\frac{d}{dt}}U, \Omega)]+g(D_{E_{i}}D_{\frac{d}{dt}}U, D_{E_{i}}\Omega)\Big\}\circ i_{t} v_{g}\nonumber\\
&=&\int_{M}\Big\{D_{E_{i}}E_{i}[g(D_{\frac{d}{dt}}U, \Omega)]-E_{i}E_{i}[g(D_{\frac{d}{dt}}U, \Omega)]-g(D_{\frac{d}{dt}}U, D_{D_{E_{i}}E_{i}}\Omega)\nonumber\\
&&+E_{i}[g(D_{\frac{d}{dt}}U, D_{E_{i}}\Omega)]+g(D_{E_{i}}D_{\frac{d}{dt}}U, D_{E_{i}}\Omega)\Big\}\circ i_{t}v_{g}\nonumber\\
&=&\int_{M}\Big\{\Delta[g(D_{\frac{d}{dt}}U, \Omega)]-g(D_{\frac{d}{dt}}U, D_{D_{E_{i}}E_{i}}\Omega)+2E_{i}[g(D_{\frac{d}{dt}}U, D_{E_{i}}\Omega)]\nonumber\\
&&-g(D_{\frac{d}{dt}}U, D_{E_{i}}D_{E_{i}}\Omega)\Big\}\circ i_{t}v_{g}\nonumber\\
&=&\int_{M}\Big\{\Delta [g(D_{\frac{d}{dt}}U, \Omega)]-g(D_{\frac{d}{dt}}U, D_{D_{E_{i}}E_{i}}\Omega)+2E_{i}[g(D_{\frac{d}{dt}}U, D_{E_{i}}\Omega)]\nonumber\\
&&-2g(D_{\frac{d}{dt}}U, D_{D_{E_{i}}E_{i}}\Omega)+2g(D_{\frac{d}{dt}}U, D_{D_{E_{i}}E_{i}}\Omega)-\nonumber\\
&&-g(D_{\frac{d}{dt}}U, D_{E_{i}}D_{E_{i}}\Omega)\Big\}\circ i_{t} v_{g}\nonumber
\end{eqnarray}
Applying the divergence Theorem for the 1-form
$\theta_{t}(\cdot)=g(\nabla_{\frac{d}{dt}}X_{t}, \nabla_{\cdot}\bar{\Delta}X_{t}), t \in
I$, we have
\begin{eqnarray}\label{Eq: derivative of ÄX}
\int_{M}g(D_{\frac{d}{dt}}\Omega, \Omega)\circ
i_{t}v_{g}&=&\int_{M}\Delta[g(\nabla_{\frac{d}{dt}}X_{t}, \bar{\Delta}X_{t})]v_{g}+2\int_{M}div \theta_{t} v_{g}+\int_{M}g(\nabla_{\frac{d}{dt}}X_{t}, \bar{\Delta}\bar{\Delta}X_{t})v_{g}\nonumber\\
&=&\int_{M}g(\nabla_{\frac{d}{dt}}X_{t}, \bar{\Delta}\bar{\Delta}X_{t})v_{g}.
\end{eqnarray}
Substituting (\ref{Eq: final derivative of Z}) and (\ref{Eq: derivative of ÄX}) in (\ref{Eq: derivative E2}), evaluating at $t=0$ and setting $V=\nabla_{\frac{d}{dt}}X_{t}|_{t=0}$, we easily obtain (\ref{Eq: variational formula}).
\end{pf}
Since the vector field $X$ is biharmonic if and only if $\frac{d}{dt}E_{2}(X_{t})|_{t=0}=0$ for all admissible variations, we get
\begin{cor}
A vector field $X$ of an $m$-dimensional Riemannian manifold $(M, g)$ is biharmonic if and only if
\begin{equation}\label{Eq: biharmonic vector field}
\bar{\Delta}\bar{\Delta}X+\sum_{i=1}^{m}[(\nabla_{e_{i}}R)(e_{i}, S(X))X+R(e_{i},
\nabla_{e_{i}}S(X))X+2R(e_{i}, S(X))\nabla_{e_{i}}X]=0,
\end{equation}
where $\{e_{i}\}_{i=1}^{m}$ is a local orthonormal frame field of $(M, g)$.
\end{cor}
\begin{rem}
Theorem 3.2 holds if $(M, g)$ is a non-compact Riemannian manifold. Indeed, if $M$ is non-compact, you must take an open subset $W$ in $M$ whose closure is compact, and take an arbitrary $V$ but, the support of $V$, namely, the closure of the set of all points in $M$ at which $V$ is not zero, is contained in $W$. Then, (\ref{Eq: variational formula}) holds in the form:
\begin{eqnarray}
\frac{d}{dt}E_{2}(X_{t})_{t=0}&=&\int_{W}\Big\{g(\bar{\Delta}\bar{\Delta}X+\sum_{i=1}^{m}\Large[(\nabla_{e_{i}}R)(e_{i}, S(X))X+R(e_{i},
\nabla_{e_{i}}S(X))X+\nonumber\\
&&+2R(e_{i}, S(X))\nabla_{e_{i}}X\Large], V)\Big\}v_{g}\nonumber
\end{eqnarray}
for such any smooth 1-parameter variation $U$. Then, you get that the term $\sum_{i=1}^{m}[(\nabla_{e_{i}}R)(e_{i}, S(X))X+R(e_{i},
\nabla_{e_{i}}S(X))X+2R(e_{i}, S(X))\nabla_{e_{i}}X]+\bar{\Delta}\bar{\Delta}X=0$ vanishes. Then, we obtain Corollary 3.3 in the case which $(M, g)$ is non-compact.
\end{rem}
\begin{rem}
In \cite{GilMedr01}, Gil-Medrano proved that a vector field $X$ in a Riemannian manifold $(M, g)$ is a harmonic vector field (i.e. critical point of the energy functional $E_{1}$, only considering variations among maps defined by vector fields) if and only if
\begin{equation*}
\bar{\Delta}X=0.
\end{equation*}
The above condition involves only the connection on the Riemannian manifold $(M, g)$. On the contrary, the corresponding critical point condition which characterizes the biharmonic vector fields involves additionally the Riemann curvature tensor of $(M, g)$.
\end{rem}
In the following Theorem, we investigate the condition under of which a vector field $X$ of a Riemannian manifold $(M, g)$ is biharmonic under the assumption that the base manifold $(M, g)$ is compact. In particular, we have
\begin{th}
Let $(M, g)$ be a compact oriented $m$-dimensional Riemannian manifold and $X\in{\mathfrak{X}}(M)$ a tangent vector field. Then, $X$ is a biharmonic vector field if and only if $X$ is parallel.
\end{th}
\begin{pf}
We assume that the vector field $X$ is a biharmonic vector field i.e. critical point of the bienergy functional $E_{2}$ restricted to the space of all vector fields of $(M, g)$. We consider the smooth 1-parameter variation $X_{t}=(1+t)X$ of $X$ ($t\in I=(-\epsilon, \epsilon), \epsilon > 0$). By using relations (\ref{Eq: helpful lemma}), (\ref{Eq: variational formula}) and the symmetries of the Riemann curvature tensor, we have
\begin{eqnarray}
0&=&\frac{d}{dt}E_{2}(X_{t})|_{t=0}=\int_{M}g(\bar{\Delta}\bar{\Delta}X, X)v_{g}+\int_{M}g\Big(\sum_{i=1}^{m}(\nabla_{e_{i}}R)(e_{i}, S(X))X, X\Big)v_{g}+\nonumber\\
&&+\int_{M}g\Big(\sum_{i=1}^{m}R(e_{i}, \nabla_{e_{i}}S(X))X, X\Big)v_{g}+2\int_{M}g\Big(\sum_{i=1}^{m}R(e_{i}, S(X))\nabla_{e_{i}}X, X\Big)v_{g}\nonumber\\
&=&\int_{M}\Delta[g(\bar{\Delta}X, X)]v_{g}+\int_{M}g(\bar{\Delta}X, \bar{\Delta}X)v_{g}+2\int_{M}div \theta_{\bar{\Delta}XX}v_{g}+\nonumber\\
&&+2\int_{M}g(S(X), S(X))v_{g}=\int_{M}\Big[g(\bar{\Delta}X, \bar{\Delta}X)+2g(S(X), S(X))\Big]v_{g},\nonumber
\end{eqnarray}
where $\{e_{i}\}_{i=1}^{m}$ a local orthonormal frame field of $(M, g)$. Furthermore, we have applied the divergence Theorem for the function $g(\bar{\Delta}X, X)$ and the 1-form $\theta_{\bar{\Delta}XX}$. Since both functions $g(\bar{\Delta}X, \bar{\Delta}X)$ and $g(S(X), S(X))$ are positive, we easily conclude that $\bar{\Delta}X=S(X)=0$ everywhere on $M$. Equivalently, $X: (M, g)\mapsto (TM, g_{S})$ is an harmonic map. By a result of Nouhaud (\cite{Nouhaud77}), we get that $X$ is parallel. Conversely, we assume that the vector field $X$ is parallel. In this case, $X : (M, g)\mapsto(TM, g_{S})$ is an harmonic map (\cite{Nouhaud77}) and, hence trivially a biharmonic map. As a consequence, $X$ is a critical point of the bienergy functional $E_{2}$ restricted to the set of all vector fields of $(M, g)$.
\end{pf}
\begin{cor}
Let $(M, g)$ be a compact oriented $m$-dimensional Riemannian manifold and $X\in{\mathfrak{X}}(M)$ a tangent vector field. Then, $X:(M, g)\mapsto (TM, g_{S})$ is a biharmonic map if and only if $X$ is parallel.
\end{cor}
\begin{rem}
In \cite{GilMedr01}, Gil-Medrano proved that if $(M, g)$ is compact, then $X$ is critical point of
the energy functional $E_{1}$ restricted to the set of all vector fields of $(M, g)$ (i.e. $X$ is a harmonic
vector field) if and only if $X$ is parallel. Theorem 3.6 shows that this result of Gil-Medrano
remains invariant in the case which the energy functional is substituted by the bienergy functional.
Furthermore, Nouhaud (\cite{Nouhaud77}) proved that if $M$ is compact, then $X : (M, g)\mapsto(TM, g_{S})$ is an
harmonic map if and only if X is parallel (see also \cite{Ishiha79}). Corollary 3.7 shows that the result of
Nouhaud remains also invariant if $X : (M, g)\mapsto (TM, g_{S})$ is biharmonic map.
\end{rem}
In the following, we give examples of non-parallel biharmonic vector fields in the case which the base manifold is non-compact. More precisely, we have
\begin{exam}
Let $(M={\mathbb{R}}^{2}, g)$ with the cartesian coordinates $(x, y)$, equipped with the standard Euclidean metric $g$. The general form of a vector field $X$ of ${\mathbb{R}}^{2}$ is
\begin{equation*}
X=f(x, y)\frac{\partial}{\partial x}+g(x, y)\frac{\partial}{\partial y},
\end{equation*}
where $f, g$ are smooth functions of ${\mathbb{R}}^{2}$ and $\frac{\partial}{\partial x}$, $\frac{\partial}{\partial y}$ are the basic vector fields of ${\mathbb{R}}^{2}$. We have
\begin{eqnarray}
\bar{\Delta}X&=&-(f_{xx}+f_{yy})\frac{\partial}{\partial x}-(g_{xx}+g_{yy})\frac{\partial}{\partial y},\nonumber\\
\bar{\Delta}\bar{\Delta}X&=&-(f_{xxxx}+2f_{xxyy}+f_{yyyy})\frac{\partial}{\partial x}-(g_{xxxx}+2g_{xxyy}+g_{yyyy})\frac{\partial}{\partial y},\nonumber
\end{eqnarray}
where $f_{x}=\frac{\partial f}{\partial x}, f_{xx}=\frac{\partial^{2}f}{\partial x^{2}}$ etc. By using (\ref{Eq: biharmonic vector field}), we deduce that the vector field $X$ is biharmonic if and only if $\bar{\Delta}\bar{\Delta}X=0$ or, equivalently,
\begin{equation}\label{Eq: equations biharmonic vector fields}
\begin{array}{lr}
\Delta^{2}f=f_{xxxx}+2f_{xxyy}+f_{yyyy}=0, & \Delta^{2}g=g_{xxxx}+2g_{xxyy}+g_{yyyy}=0.
\end{array}
\end{equation}
Particular solutions of the system (\ref{Eq: equations biharmonic vector fields}) are
\begin{eqnarray}\label{Eq: solutions biharmonic vector fields}
f(x, y)&=&\{(A+Cx)\cosh \beta x+(B+Dx)\sinh \beta x\}\{a\cos\beta y+b\sin \beta y\},\nonumber\\
g(x, y)&=&\{(A+Cx)\cos\beta x+(B+Dx)\sin \beta x\}\{a \cosh \beta y+b \sinh \beta y\},\nonumber\\
\end{eqnarray}
where $A, B, C, D, a, b$ and $\beta$ are arbitrary real constants.
\end{exam}
\begin{rem}
The equation $\Delta^{2}f=0$ is called \emph{biharmonic equation} and the corresponding solutions are called \emph{biharmonic functions}. The biharmonic equation is encountered in plane problems of elasticity ($f$ is the Airy stress function). It is also used to describe slow flows of viscous incompressible fluids ($f$ is the stream function). For more information, we refer to the book \cite{Poly02}. Considering the function $f(x, y)$ of the system (\ref{Eq: equations biharmonic vector fields}) given in (\ref{Eq: solutions biharmonic vector fields}), we calculate
\begin{equation*}
f_{xx}+f_{yy}=2\beta\{D\cosh \beta x + C\sinh \beta x\}\{a\cos \beta y+b\sin \beta y\}.
\end{equation*}
Considering the case $\beta=C=D=a=1$ and $b=0$, we obtain $f_{xx}+f_{yy}=2e^{x}\cos y\neq 0$ for all $(x, y)\in {\mathbb{R}}^{2}-\{(t, 2\kappa\pi+\frac{\pi}{2})/t\in{\mathbb{R}},\kappa \in {\mathbb{Z}}\}$. As a conclusion, for this particular choice of the parameters $\beta, C, D, a$ and $b$, the corresponding solution $f(x, y)$ of (\ref{Eq: equations biharmonic vector fields}) is a biharmonic function but not harmonic function. Similarly, for the solution $g(x, y)$ of the system (\ref{Eq: equations biharmonic vector fields}) given in (\ref{Eq: solutions biharmonic vector fields}), we have
\begin{equation*}
g_{xx}+g_{yy}=2\beta\{D\cos\beta x-C\sin \beta x\}\{a\cosh \beta y+b\sinh \beta y\}.
\end{equation*}
Considering the case $a=b=\beta=D=1$ and $C=0$, we obtain $g_{xx}+g_{yy}=2e^{y}\cos x\neq 0$ for all $(x, y)\in {\mathbb{R}}^{2}-\{(2\kappa\pi+\frac{\pi}{2}, t)/\kappa\in{\mathbb{Z}}, t\in{\mathbb{R}}\}$. Summarizing, we yield that the vector fields $X$ of ${\mathbb{R}}^{2}$
\begin{equation*}
X=\{(A+x)\cosh x+(B+x)\sinh x\}\cos y\frac{\partial}{\partial x}+e^{y}\{A\cos x+(B+x)\sin x\}\frac{\partial}{\partial y},
\end{equation*}
where $A, B$ are arbitrary real constants, are non-parallel biharmonic vector fields but non harmonic (see Remark 3.5).
\end{rem}
\begin{rem}
Given $u(x, y)$ and $v(x, y)$ solutions of the Laplace equation ($\Delta w=0$), we have the following various representations of the general solution of the biharmonic equation $\Delta^{2}w=0$ (\cite[p. 516]{Poly02}):
\begin{eqnarray}
w(x, y)&=&x u(x, y)+v(x, y)\nonumber\\
w(x, y)&=&y u(x, y)+v(x, y)\nonumber\\
w(x, y)&=&(x^{2}+y^{2})u(x, y)+v(x, y).\nonumber
\end{eqnarray}
Equivalently, we have the following complex form of the presentation of the general solution:
\begin{equation*}
w(x, y)=Re[\bar{z}f(z)+g(z)],
\end{equation*}
where $f(z)$ and $g(z)$ are arbitrary analytic functions of the variable $z=x+i y, \bar{z}=x-i y$ (see also \cite{AroCreeLipkin83}). As a consequence, the general representation formula of a biharmonic vector field $X$ of $({\mathbb{R}}^{2}, g)$ is
\begin{equation*}
X=(x u(x, y)+v(x, y))\frac{\partial}{\partial x}+(y u(x, y)+v(x, y))\frac{\partial}{\partial y},
\end{equation*}
where $u(x, y)$ and $v(x, y)$ are arbitrary harmonic functions.
\end{rem}
\begin{exam}
We consider the Heisenberg group $Nil_{3}$. This is defined to be the group consisting of all real $3\times 3$ upper triangular matrices of the form
\begin{equation*}
\begin{array}{c}
A
\end{array}
=\left (
\begin{array}{cccc}
1 &x &y\\
0 &1 &z\\
0 & 0 & 1
\end{array}\right )
\end{equation*}
endowed with the left-invariant metric given by $(dx)^{2}+(dy-xdz)^{2}+(dz)^{2}$. We may thus identify $Nil_{3}$ with ${\mathbb{R}}^{3}$ endowed with this metric. The left-invariant vector fields
\begin{equation*}
\begin{array}{lcr}
e_{1}=\frac{\partial}{\partial x}, &e_{2}=\frac{\partial}{\partial y}, &e_{3}=\frac{\partial}{\partial z}+x\frac{\partial}{\partial y},
\end{array}
\end{equation*}
constitute an orthonormal basis of the Lie algebra $\mathfrak{g}$ of $Nil_{3}$. The corresponding Levi Civita connection is determined by
\begin{equation}\label{Eq: connection Heisenberg group}
\begin{array}{lcr}
\nabla_{e_{1}}e_{2}=\nabla_{e_{2}}e_{1}=-\frac{1}{2}e_{3}, & \nabla_{e_{1}}e_{3}=-\nabla_{e_{3}}e_{1}=\frac{1}{2}e_{2}, & \nabla_{e_{2}}e_{3}=\nabla_{e_{3}}e_{2}=\frac{1}{2}e_{1},
\end{array}
\end{equation}
where the remaining covariant derivatives of the basic vectors vanish. By using relations (\ref{Eq: connection Heisenberg group}), we get
\begin{equation*}
R(e_{2}, e_{1})e_{3}=R(e_{2}, e_{3})e_{1}=R(e_{3}, e_{1})e_{2}=0,
\end{equation*}
\begin{equation}\label{Eq: Heisenberg group 1}
\begin{array}{lcr}
\bar{\Delta}e_{1}=-\nabla_{e_{1}}\nabla_{e_{1}}e_{1}-\nabla_{e_{2}}\nabla_{e_{2}}e_{1}-\nabla_{e_{3}}\nabla_{e_{3}}e_{1}=\frac{1}{2}e_{1}, & \bar{\Delta}e_{2}=\frac{1}{2}e_{2}, & \bar{\Delta}e_{3}=\frac{1}{2}e_{3}.
\end{array}
\end{equation}
We consider the vector field $X=f(x)e_{1}$, where $f(x)$ is a smooth function of ${\mathbb{R}}$ depending of the variable $x$. Combining relations (\ref{Eq: tensorial character S(X)}), (\ref{Eq: rough Laplacian fX}), (\ref{Eq: connection Heisenberg group}) and (\ref{Eq: Heisenberg group 1}), we obtain
\begin{eqnarray}
\bar{\Delta}X&=&\Delta f e_{1}+f\bar{\Delta}e_{1}-2f^{'}\nabla_{e_{1}}e_{1}=(\frac{1}{2}f-f^{''})e_{1},\nonumber\\
\bar{\Delta}\bar{\Delta}X&=&\Delta(\frac{1}{2}f-f^{''})e_{1}+(\frac{1}{2}f-f^{''})\bar{\Delta}e_{1}=(f^{''''}-f^{''}+\frac{1}{4}f)e_{1},
\nonumber\\
S(X)&=&f^{2}S(e_{1})=f^{2}(-\frac{1}{2}R(e_{3}, e_{1})e_{2}-\frac{1}{2}R(e_{2}, e_{1})e_{3})=0,\nonumber
\end{eqnarray}
where $f^{'}=\frac{df}{d x}, f^{''}=\frac{d^{2}f}{d x^{2}}$ etc.By using (\ref{Eq: biharmonic vector field}), we deduce that the vector field $X$ is biharmonic if and only if $\bar{\Delta}\bar{\Delta}X=0$ or, equivalently,
\begin{equation}\label{Eq: biharmonic Heisenberg group}
f^{''''}-f^{''}+\frac{1}{4}f=0.
\end{equation}
The general solution of the homogeneous fourth order equation (\ref{Eq: biharmonic Heisenberg group}) is
\begin{equation*}
f(x)=c_{1}e^{\frac{x}{\sqrt{2}}}+c_{2}xe^{\frac{x}{\sqrt{2}}}+c_{3}e^{-\frac{x}{\sqrt{2}}}+c_{4}xe^{-\frac{x}{\sqrt{2}}},
\end{equation*}
where $c_{1}, c_{2}, c_{3}$ and $c_{4}$ are real constants. We easily conclude that the vector fields $X=[c_{1}e^{\frac{x}{\sqrt{2}}}+c_{3}e^{-\frac{x}{\sqrt{2}}}]e_{1}, c_{1}, c_{3}\in {\mathbb{R}}$ are harmonic and, also, biharmonic vector fields.  On the contrary, the vector fields $X=x(c_{2}e^{\frac{x}{\sqrt{2}}}+c_{4}e^{-\frac{x}{\sqrt{2}}})e_{1}$ are biharmonic but non harmonic vector fields of $Nil_{3}$. Following the same procedure, the vector fields $X=z(c_{2}e^{\frac{z}{\sqrt{2}}}+c_{4}e^{-\frac{z}{\sqrt{2}}})e_{3}, c_{2}, c_{4}\in {\mathbb{R}}$ are also biharmonic but non harmonic vector fields of $Nil_{3}$.
\end{exam}
In the following example, we give examples of vector fields which are harmonic vector fields but not biharmonic vector fields. More precisely, we have
\begin{exam}
We consider the hyperbolic space $(H^{n}, g)$ of constant sectional curvature $-c^{2}<0$
\begin{eqnarray}
H^{n}&=&\Large\{z=(x_{1},\ldots, x_{n-1}, y)\in {\mathbb{R}}^{n}: y>0\Large\},\nonumber\\
g&=&\frac{1}{(cy)^{2}}\Big (\sum_{i=1}^{n-1}(dx_{i})^{2}+(dy)^{2}\Big).\nonumber
\end{eqnarray}
The vector fields
\begin{equation*}
\begin{array}{lcr}
V=cy \frac{\partial}{\partial y}, & E_{i}=cy\frac{\partial}{\partial x_{i}}, & 1\leq i \leq n-1,
\end{array}
\end{equation*}
constitute an orthonormal frame field of $(H^{n}, g)$ and
\begin{equation*}
\begin{array}{lcr}
[V, E_{i}]=cE_{i}, & [E_{i}, E_{j}]=0, & 1 \leq i, j \leq n-1.
\end{array}
\end{equation*}
Let $\nabla$ be the Levi-Civita connection of $(H^{n}, g)$. It follows that (\cite[p. 139]{DragomirPerro11})\\
\begin{equation}\label{Eq: connection H2}
\begin{array}{llcr}
\nabla_{E_{i}}V=-cE_{i}, &\nabla_{V}{V}=0, &\nabla_{V}E_{i}=0, &\nabla_{E_{i}}E_{j}=c\delta_{ij}V,
\end{array}
\end{equation}
for any $1 \leq i, j \leq n-1$. Then (\cite[p. 139]{DragomirPerro11})\\
\begin{equation}\label{Eq: connection H2 I}
\begin{array}{ll}
\bar{\Delta}V=(n-1)c^{2}V, & S(V)=-c^{2}(\nabla_{V}{V}-div (V) V)=-c^{3}(n-1)V, \\
div (V)=-(n-1)c.
\end{array}
\end{equation}
We consider the vector field $X=f(y)V$, where $f(y)$ is a smooth function of ${\mathbb{R}}^{+}$ depending on the variable $y$ ($y>0$). By using relations (\ref{Eq: tensorial character S(X)}), (\ref{Eq: rough Laplacian fX}), (\ref{Eq: connection H2}) and (\ref{Eq: connection H2 I}), we get
\begin{eqnarray}\label{Eq: biharmonicity H2 I}
\bar{\Delta}X&=&\Delta f V + f\bar{\Delta}V-2V(f)\nabla_{V}V=[(n-1)c V(f)-V V(f)]V + (n-1)c^{2}f V\nonumber\\
&=&c^{2}\Big((n-2)yf^{'}-y^{2}f^{''}+(n-1)f\Big)V,\nonumber\\
\bar{\Delta}\bar{\Delta}X&=&c^{2}\Big (\Delta [(n-2)yf^{'}-y^{2}f^{''}+(n-1)f]V+[(n-2)yf^{'}-y^{2}f^{''}+(n-1)f]\bar{\Delta}V\Big)\nonumber\\
&=&c^{4}\Big (y^{4}f^{''''}+(8-2n)y^{3}f^{'''}+(n-8)(n-2)y^{2}f^{''}+(n-2)(3n-4)yf^{'}+\nonumber\\
&&+(n-1)^{2}f\Big)V, \hspace{0.10cm} S(X)=f^{2}S(V)=-c^{3}(n-1)f^{2}V,
\end{eqnarray}
where $f^{'}=\frac{df}{dy}, f^{''}=\frac{d^{2}f}{dy^{2}}$ etc. By using relations (\ref{Eq: connection H2}), (\ref{Eq: biharmonicity H2 I}) and the fact that $(H^{n}, g)$ is a space of constant curvature $-c^{2}$, we have
\begin{equation*}
\begin{array}{lcr}
\nabla_{E_{i}}X=-c f E_{i}, & \nabla_{V}S(X)=-c^{3}(n-1) V(f^{2})V, & \nabla_{E_{i}}S(X)=c^{4}(n-1)f^{2}E_{i}, \\
\end{array}
\end{equation*}
\begin{equation}\label{Eq: biharmonicity H2 II}
\begin{array}{lcr}
R(V, \nabla_{V}S(X))V=0, & R(E_{i}, \nabla_{E_{i}}S(X))X=0, & R(V, S(X))\nabla_{V}X=0, \\
\end{array}
\end{equation}
\begin{equation*}
\begin{array}{lr}
\sum_{i=1}^{n-1}R(E_{i}, S(X))\nabla_{E_{i}}X=c^{6}(n-1)^{2}f^{3}V, & \nabla R = 0.
\end{array}
\end{equation*}
for any $1\leq i \leq n-1$. By using the first relation of (\ref{Eq: biharmonicity H2 I}), we easily conclude that the vector field $X=f(y)V$ is a harmonic vector field (equivalently, $\bar{\Delta}X=0$) if and only if
\begin{equation}\label{Eq: harmonic H2}
y^{2}f^{''}-(n-2)yf^{'}-(n-1)f=0.
\end{equation}
The ODE (\ref{Eq: harmonic H2}) is a second order Euler's equation and its general solution is
\begin{equation*}
f(y)=c_{1}y^{\frac{n-1-\sqrt{(n-1)(n+3)}}{2}}+c_{2}y^{\frac{n-1+\sqrt{(n-1)(n+3)}}{2}},\hspace{0.05cm} y>0
\end{equation*}
where $c_{1}$ and $c_{2}$ are real constants. On the other hand, combining relations (\ref{Eq: biharmonic vector field}), (\ref{Eq: biharmonicity H2 I}) and (\ref{Eq: biharmonicity H2 II}), we easily conclude that the vector field $X=f(y)V$ is biharmonic if and only if
\begin{eqnarray}\label{Eq: biharmonic equation H2}
y^{4}f^{''''}+(8-2n)y^{3}f^{'''}+(n-8)(n-2)y^{2}f^{''}+(n-2)(3n-4)yf^{'}+\nonumber\\
+(n-1)^{2}f=-2c^{2}(n-1)^{2}f^{3}.\nonumber\\
\end{eqnarray}
Clearly, the vector fields $X=y^{\frac{n-1-\sqrt{(n-1)(n+3)}}{2}}V$ and $X=y^{\frac{n-1+\sqrt{(n-1)(n+3)}}{2}}V$ are harmonic vector fields but not biharmonic vector fields.
\end{exam}
\begin{rem}
The fourth order ODE (\ref{Eq: biharmonic equation H2}) has a particular interest. Applying the substitution $y=e^{t}$, equation (\ref{Eq: biharmonic equation H2}) is transformed to
\begin{eqnarray}\label{Eq: biharmonic equation H2 new}
v^{''''}+(2-2n)v^{'''}+(n-1)(n-3)v^{''}+2(n-1)^{2}v^{'}+(n-1)^{2}v+\nonumber\\
+2c^{2}(n-1)^{2}v^{3}=0,\nonumber\\
\end{eqnarray}
where $v(t)=f(y)$ and $t\in {\mathbb{R}}$. The ODE (\ref{Eq: biharmonic equation H2 new}) has no global solution, since the term $v^{3}$ goes up very quickly at infinity.
However, the ODE (\ref{Eq: biharmonic equation H2 new}) has a local solution.
\end{rem}
\begin{rem}
Examples 3.9, 3.12 and 3.13 show that there exist biharmonic vector fields which are not harmonic, and conversely. As a consequence, the notions ''biharmonic  vector field'' and ''harmonic vector field'' are independent in the sense that it does not exist an immediate relation between them.
\end{rem}
\subsection*{Acknowledgment}
The authors acknowledge Professor Hisashi Naito for his suggestion about the non existence of global solution of the ODE (\ref{Eq: biharmonic equation H2}).

\vskip0.6cm\par
Pafsania, Isthmia Korinthias,
\par Korinthos, GR-20010, Greece.
\par
{\it E-mail address}: mmarkellos@@hotmail.gr
\vskip0.8cm\par
Institute for International Education,
\par Tohoku University,
Kawauchi 41, Sendai, 980-8576, Japan.
\par
{\it E-mail address}: urakawa@@math.is.tohoku.ac.jp

\begin{thebibliography}{99}

\bibitem{AroCreeLipkin83}
N.~Aronszajn, T.~M.~Creese and L.~J.~Lipkin,
\newblock{\em Polyharmonic functions},
\newblock Oxford, 1983.

\bibitem{Baird03}
P.~Baird and J.~C.~Wood,
\newblock{\em Harmonic morphisms between Riemannian manifolds},
\newblock Oxford University Press, {\bf{29}}, London Math. Soc.
Monogr. (N.S.), 2003.

\bibitem{CaMoOn01}
R.~Caddeo, S.~Montaldo, and C.~Oniciuc,
\newblock {\em Biharmonic submanifolds of ${\mathbb{S}}^{3}$},
\newblock Intern. J. Math, {\bf{12}} (2001), 867--876.

\bibitem{CaMoOn02}
R.~Caddeo, S.~Montaldo, and C.~Oniciuc,
\newblock{\em Biharmonic submanifolds in spheres},
\newblock Israel J.Math., {\bf{130}} (2002), 109--123.

\bibitem{DjaaElheOuak12}
Mustapha Djaa, Hichem Elhendi and Seddik Ouakkas,
\newblock{\em On the biharmonic vector fields},
\newblock Turk. J. Math., {\bf{36}} (2012), 463--474.

\bibitem{DragomirPerro11}
S. Dragomir and D. Perrone,
\newblock {\em Harmonic Vector Fields: Variational Principles and
Differential Geometry}, Elsevier 2011.

\bibitem{EellsSampson64}
J.~Eells and J.~H.~Sampson,
\newblock {\em Harmonic mappings of Riemannian manifolds},
Amer. J. Math., {\bf{86}}(1) (1964), 109--160.

\bibitem{GilMedr01}
O.~Gil-Medrano,
\newblock{\em Relationship between volume and energy of
unit vector fields},
Differential Geom.  Appl., {\bf{15}} (2001), 137--152.

\bibitem{GudmuKappos02}
S.~Gudmundsson and E.~Kappos,
\newblock {\em On the Geometry of Tangent Bundles},
Expo. Math., {\bf{20}} (2002), 1--41.

\bibitem{HiKaWo01}
A.~Higuchi, B.~S.~Kay and C.~M.~Wood,
\newblock{\em The energy of unit vector fields on the 3-sphere},
J.Geom.Phys., {\bf{37}}(1-2) (2001), 137--155.

\bibitem{Ishiha79}
T.~Ishihara,
\newblock{\em Harmonic sections of tangent bundles},
J. Math. Tokushima Univ., {\bf{13}} (1979), 23--27.

\bibitem{Jiang08}
G.~Jiang,
\newblock{\em 2-harmonic maps and their first and second variational
formulas}, Translated into English by Hajime Urakawa,
Note Mat., {\bf{28}} suppl. n. 1 (2008), 209--232.

\bibitem{MontaOni06}
S.~Montaldo and C.~Oniciuc,
\newblock{\em A short survey on biharmonic maps between Riemannian
manifolds},
Rev. Un. Mat. Argentina, {\bf{47}}(2) (2006), 1--22.

\bibitem{Nouhaud77}
O.~Nouhaud,
\newblock{\em Applications harmoniques d' une vari\'{e}t\'{e} Riemannienne dans son fibr\'{e}
tangent},
C.R.Acad.Sci.Paris, {\bf{284}} (1977), 815--818.

\bibitem{Poly02}
A.~D.~Polyanin,
\newblock {\em Handbook of Linear Partial Differential Equations for Engineers and Scientists},
\newblock Chapman and Hall/CRC, 2002.

\bibitem{Urakawa}
H.~Urakawa,
\newblock{\em Calculus of Variations and Harmonic Maps},
\newblock Amer. Math. Soc., Providence {\bf{132}}, Transl. Math.
Monograph, 1993.

\end{thebibliography}
\end{document}